\documentclass[review,5p,times,twocolumn,round,authoryear,square]{elsarticle}

\usepackage{amssymb,amsmath}
\usepackage{graphics,graphicx}
\usepackage{url,theorem}
\usepackage{booktabs}
\usepackage{xcolor}

\usepackage{epstopdf}

\journal{Chemical Engineering Science}

\newcommand{\torol}[1]{}
\newcommand{\irrevreak}[3]{#1{\ }{\overset{#2}{\longrightarrow}}{\ }#3}
\newcommand{\revreak}[4]{#1{\ }{\overset{#2}{\underset{#3}{\rightleftharpoons}}#4}}

\theorembodyfont{\upshape}

\newcommand{\ikde}{induced kinetic differential equation}
\newcommand{\ivp}{initial value problem}
\newcommand{\Mma}{\textit{Mathematica}}
\newcommand{\ode}{ordinary differential equation}

\newcommand{\Vg}{Volpert graph}
\newcommand{\AbsoluteRobustness} {\ifmmode {\textbf{\texttt{AbsoluteRobustness}}} \else {\bf\tt AbsoluteRobustness}\fi}
\newcommand{\AtomMatrix}         {\ifmmode {\textbf{\texttt{AtomMatrix}}} \else {\bf\tt AtomMatrix}\fi}
\newcommand{\Bold}               {\ifmmode {\textbf{\texttt{Bold}}} \else {\bf\tt Bold}\fi}
\newcommand{\Break}              {\ifmmode {\textbf{\texttt{Break}}} \else {\bf\tt Break}\fi}
\newcommand{\Compile}            {\ifmmode {\textbf{\texttt{Compile}}} \else {\bf\tt Compile}\fi}
\newcommand{\Continue}           {\ifmmode {\textbf{\texttt{Continue}}} \else {\bf\tt Continue}\fi}
\newcommand{\CUDALink}           {\ifmmode {\textbf{\texttt{CUDALink}}} \else {\bf\tt CUDALink}\fi}
\newcommand{\DirectedEdges}      {\ifmmode {\textbf{\texttt{DirectedEdges}}} \else {\bf\tt DirectedEdges}\fi}
\newcommand{\Do}                 {\ifmmode {\textbf{\texttt{Do}}} \else {\bf\tt Do}\fi}
\newcommand{\EdgeLabeling}       {\ifmmode {\textbf{\texttt{EdgeLabeling}}} \else {\bf\tt EdgeLabeling}\fi}
\newcommand{\ElementaryReactions}{\ifmmode {\textbf{\texttt{ElementaryReactions}}} \else {\bf\tt ElementaryReactions}\fi}
\newcommand{\ExponentialDistribution}{\ifmmode {\textbf{\texttt{ExponentialDistribution}}} \else {\bf\tt ExponentialDistribution}\fi}
\newcommand{\False}              {\ifmmode {\textbf{\texttt{False}}} \else {\bf\tt False}\fi}
\newcommand{\FindFit}              {\ifmmode {\textbf{\texttt{FindFit}}} \else {\bf\tt FindFit}\fi}
\newcommand{\For}                {\ifmmode {\textbf{\texttt{For}}} \else {\bf\tt For}\fi}
\newcommand{\GetProblem}         {\ifmmode {\textbf{\texttt{GetProblem}}} \else {\bf\tt GetProblem}\fi}
\newcommand{\Goto}               {\ifmmode {\textbf{\texttt{Goto}}} \else {\bf\tt Goto}\fi}
\newcommand{\ImageSize}          {\ifmmode {\textbf{\texttt{ImageSize}}} \else {\bf\tt ImageSize}\fi}
\newcommand{\Label}              {\ifmmode {\textbf{\texttt{Label}}} \else {\bf\tt Label}\fi}
\newcommand{\Manipulate}         {\ifmmode {\textbf{\texttt{Manipulate}}} \else {\bf\tt Manipulate}\fi}
\newcommand{\OpenCLLink}         {\ifmmode {\textbf{\texttt{OpenCLLink}}} \else {\bf\tt OpenCLLink}\fi}
\newcommand{\Parallelize}        {\ifmmode {\textbf{\texttt{Parallelize}}} \else {\bf\tt Parallelize}\fi}
\newcommand{\ParallelMap}        {\ifmmode {\textbf{\texttt{ParallelMap}}} \else {\bf\tt ParallelMap}\fi}
\newcommand{\PlotFunction}       {\ifmmode {\textbf{\texttt{PlotFunction}}} \else {\bf\tt PlotFunction}\fi}
\newcommand{\PlotLabel}          {\ifmmode {\textbf{\texttt{PlotLabel}}} \else {\bf\tt PlotLabel}\fi}
\newcommand{\PoissonDistribution}{\ifmmode {\textbf{\texttt{PoissonDistribution}}} \else {\bf\tt PoissonDistribution}\fi}
\newcommand{\RandomChoice}       {\ifmmode {\textbf{\texttt{RandomChoice}}} \else {\bf\tt RandomChoice}\fi}
\newcommand{\RandomVariate}      {\ifmmode {\textbf{\texttt{RandomVariate}}} \else {\bf\tt RandomVariate}\fi}
\newcommand{\ReactionKinetics}   {\ifmmode {\textbf{\texttt{ReactionKinetics}}} \else {\bf\tt ReactionKinetics}\fi}
\newcommand{\RegularExpression}  {\ifmmode {\textbf{\texttt{RegularExpression}}} \else {\bf\tt RegularExpression}\fi}
\newcommand{\Return}             {\ifmmode {\textbf{\texttt{Return}}} \else {\bf\tt Return}\fi}
\newcommand{\ShowVolpertGraph}   {\ifmmode {\textbf{\texttt{ShowVolpertGraph}}} \else {\bf\tt ShowVolpertGraph}\fi}
\newcommand{\Style}              {\ifmmode {\textbf{\texttt{Style}}} \else {\bf\tt Style}\fi}
\newcommand{\Times}              {\ifmmode {\textbf{\texttt{Times}}} \else {\bf\tt Times}\fi}
\newcommand{\True}               {\ifmmode {\textbf{\texttt{True}}} \else {\bf\tt True}\fi}
\newcommand{\VertexLabeling}     {\ifmmode {\textbf{\texttt{VertexLabeling}}} \else {\bf\tt VertexLabeling}\fi}
\newcommand{\While}              {\ifmmode {\textbf{\texttt{While}}} \else {\bf\tt While}\fi}

\newcommand{\RM}{\mathbb{R}^M}

\newcommand{\alphab}{\mbox{\boldmath$\alpha$}}
\newcommand{\betab}{\mbox{\boldmath$\beta$}}
\newcommand{\gammab}{\mbox{\boldmath$\gamma$}}

\newcommand{\nulb}{\ifmmode \mathbf{0}\else \textbf{0}\fi}
\newcommand{\oneb}{\ifmmode \mathbf{1}\else \textbf{1}\fi}
\newcommand{\twob}{\ifmmode \mathbf{2}\else \textbf{2}\fi}

\newcommand{\ab}{\ifmmode \mathbf{a}\else \textbf{a}\fi}
\newcommand{\Ab}{\ifmmode \mathbf{A}\else \textbf{A}\fi}
\newcommand{\bb}{\ifmmode \mathbf{b}\else \textbf{b}\fi}
\newcommand{\Bb}{\ifmmode \mathbf{B}\else \textbf{B}\fi}
\newcommand{\cb}{\ifmmode \mathbf{c}\else \textbf{c}\fi}
\newcommand{\Cb}{\ifmmode \mathbf{C}\else \textbf{C}\fi}
\newcommand{\Db}{\ifmmode \mathbf{D}\else \textbf{D}\fi}
\newcommand{\fb}{\ifmmode \mathbf{f}\else \textbf{f}\fi}
\newcommand{\Fb}{\ifmmode \mathbf{F}\else \textbf{F}\fi}
\newcommand{\gb}{\ifmmode \mathbf{g}\else \textbf{g}\fi}
\newcommand{\hb}{\ifmmode \mathbf{h}\else \textbf{h}\fi}
\newcommand{\kb}{\ifmmode \mathbf{k}\else \textbf{k}\fi}
\newcommand{\Kb}{\ifmmode \mathbf{K}\else \textbf{K}\fi}
\newcommand{\Mb}{\ifmmode \mathbf{M}\else \textbf{M}\fi}
\newcommand{\nb}{\ifmmode \mathbf{n}\else \textbf{n}\fi}
\newcommand{\pb}{\ifmmode \mathbf{p}\else \textbf{p}\fi}
\newcommand{\Pb}{\ifmmode \mathbf{P}\else \textbf{P}\fi}
\newcommand{\qb}{\ifmmode \mathbf{q}\else \textbf{q}\fi}
\newcommand{\rb}{\ifmmode \mathbf{r}\else \textbf{r}\fi}
\newcommand{\sbold}{\ifmmode \mathbf{s}\else \textbf{s}\fi}
\newcommand{\vb}{\ifmmode \mathbf{v}\else \textbf{v}\fi}
\newcommand{\wb}{\ifmmode \mathbf{w}\else \textbf{w}\fi}
\newcommand{\xb}{\ifmmode \mathbf{x}\else \textbf{x}\fi}
\newcommand{\Xb}{\ifmmode \mathbf{X}\else \textbf{X}\fi}
\newcommand{\yb}{\ifmmode \mathbf{y}\else \textbf{y}\fi}
\newcommand{\zb}{\ifmmode \mathbf{z}\else \textbf{z}\fi}

\newcommand{\dr}{\mathrm{d}}%
\begin{document}

\begin{frontmatter}
\title{
{\ReactionKinetics}---A {\it Mathematica} package with applications II.\\
Computational problems when building a reaction kinetics package}
\author[bme]{A. L. Nagy\fnref{fn1,fn2}}
\ead{nagyal@math.bme.hu}
\author[nw]{D. Papp\fnref{fn3}}
\ead{dpapp@iems.northwestern.edu}
\author[bme,elte]{J. T\'oth\corref{cor1}\fnref{fn1,fn2}}
\ead{jtoth@math.bme.hu}

\address[bme]{Department of Analysis,
Budapest University of Technology and Economics,
Egry J. u. 1., Budapest, Hungary, H-1111}

\address[elte]{Laboratory for Chemical Kinetics of the
Institute of Chemistry,
E\"otv\"os Lor\'and University,
Budapest, Hungary}

\address[nw]{Department of Industrial Engineering and Management Sciences,\\
Northwestern University,
2145 Sheridan Road, Room C210, Evanston, IL 60208}

\fntext[fn1]{Partially supported by the Hungarian National Scientific Foundation, No. 84060.}
\fntext[fn2]{This work is connected to the scientific program of the "Development of quality-oriented and harmonized R+D+I strategy and functional model at BME" project. This project is supported by the New Sz\'echenyi Plan (Project ID: T\'AMOP-4.2.1/B-09/1/KMR-2010-0002).}
\fntext[fn3]{Partially supported by the COST Action CM901: Detailed Chemical Kinetic Models for Cleaner Combustion.}

\cortext[cor1]{Corresponding author}
\begin{abstract}
Treating a realistic problem in any field of reaction kinetics raises a series of problems:
we review these illustrated with examples
using \ReactionKinetics, a \Mma\ based package.
\end{abstract}


\end{frontmatter}

\section{Introduction}
In Part I of our paper \cite{tothnagypapp} we have formulated the requirements for
a reaction kinetics package to be useful for
a wide circle of users.

In the present Part II we enumerate the major problems arising when writing and using such a package.
It turns out that the solution of some problems is far from being trivial.
\section{Parsing}
As we mentioned in the first part of our paper, at the very beginning of application of computers to chemical kinetics the problem arouse how to construct the \ikde\ of a reaction without making to much errors.
The best thing is if the chemist provides the reaction steps, and the program creates the \ikde s automatically, if the kinetics is supposed to be of the mass action type.
If it is not then reaction rates should also be provided by the user.

As our program is based on a modern mathematical program package (often called by the misnomer
\emph{computer algebra system}) there is no limitation as to the size of the reaction to be handled,
see the Lendvay example in our previous paper.
Confer this with the quite typical restriction "\ldots has the capabilities for handling up to $590$ differential equations."

It may be instructive to cite the code building the right hand side of the \ikde\
of the reaction (note that reversibility is not assumed)
\[
\irrevreak{\sum_{m=1}^M\alpha(m,r)X(m)}{k_r}{\sum_{m=1}^M\beta(m,r)X(m)}\quad(r=1,2,\dots,R)
\]
given the matrices ($\alphab,\betab$) of molecularities, the vector of reaction rate coefficients ($\kb=\left(k_r\right)_{r=1}^R$) and the concentration vector ($\cb$) the right-hand side of the kinetic ODE gets an exceedingly transparent form (in coding as well):
\[
(\betab-\alphab).(\kb{\ }\Times\mathrm{@@}{\ }\cb^{\alphab}),
\]
where notice that both dot and componentwise products are used. 
The \emph{reaction step vector} $\gammab(\cdot,r)$, where $\gammab=\betab-\alphab$ expresses as the effect of the $r^{\mathrm{th}}$ reaction step.

\section{Combinatorics}
In reaction kinetics graphs of many types are useful and used, though it is not quite obvious how to visualize them to be the most informative for users, especially concerning to large mechanisms. Hence we accept the representations offered by \Mma.

\begin{center}
\begin{figure}[h!]\label{figglyc}
\includegraphics[height=7cm,width=10cm]{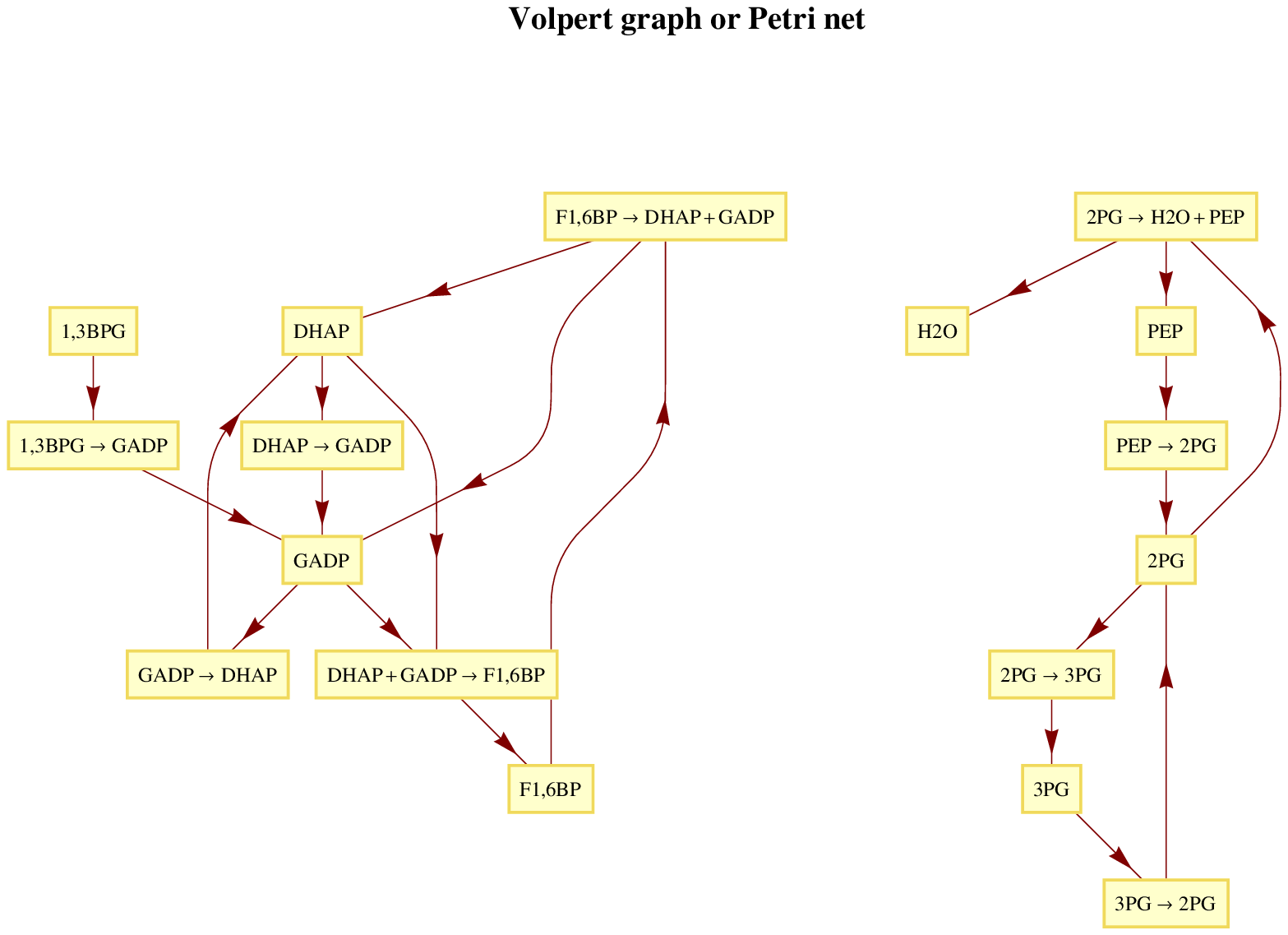}
\includegraphics[height=15cm,width=10cm]{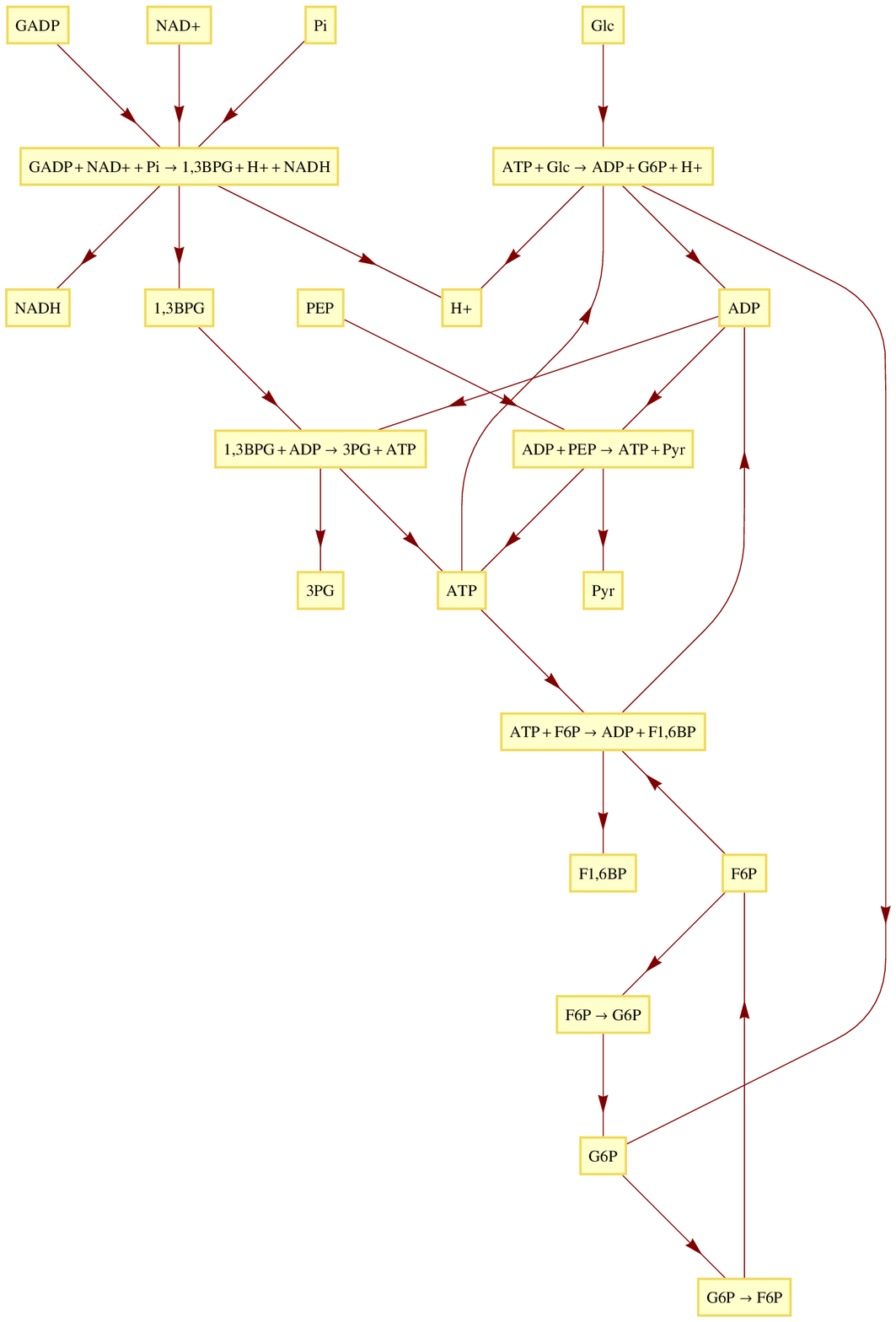}
\caption{The \Vg\ or Petri net of a glycolysis model}
\end{figure}
\end{center}

The \emph{Volpert graph} is one of the widely applied graphs in reaction kinetics which proved to be a good tool to investigate problems emerging e.g. in the theory of kinetic \ode s. 
It is a directed bipartite graph with species and reaction steps as its vertices and with $\alpha(m,r)$ edges from the vertex representing species $X(m)$ into the vertex representing reaction step $r$ and  with $\beta(m,r)$ edges from the vertex representing reaction step $r,$ into the vertex representing species $X(m)$ \citep{erditoth}.
At the same time this is the graph one can see in textbooks on biochemistry, in metabolic maps or at MaCKiE Workshops.
Without any additional advantage they can also be called \emph{Petri net}s \citep{martinezsilva}.
A model of glycolysis is a built-in reaction of our program package therefore having got it 
we can easily draw its Volpert graph which turns out to be 
not so much different from the one used in textbooks on biochemistry (see Figure \ref{figglyc}).
\vspace{-1.5cm}

\section{Linear (Diophantine) equations}
Complex reactions occur via pathways of simple ``elementary'' reaction steps. It is a common problem that the overall (or global) reaction is measured, and one would like to reconstruct the underlying network of simple reaction steps. By simple here one may mean reaction steps of order not higher than two.
Let us consider an example \citep{kovacsvizvaririedeltoth}.
We start form the overall description of the permanganate/oxalate reaction
\[
\revreak{2\mathrm{MnO}^{4-}+6\mathrm{H}^{+}+5\mathrm{C}_2\mathrm{H}_2\mathrm{O}_4}{}{}
{2\mathrm{Mn}^{2+}+8\mathrm{H}_{2}\mathrm{O}+10\mathrm{CO}_2} .
\]
The following steps lead to its the decompositions to elementary steps:
\begin{enumerate}
\item
determine the combinatorially possible species and select from those the chemically acceptable ones,
\item
determine the combinatorially possible reaction steps which are simple in the above sense and select from those the chemically acceptable ones,
\item
find those representations of the given overall reaction which are chemically acceptable.
\end{enumerate}

As we have seen in Part I of our paper, the second and third steps reduce to the solution of linear Diophantine equations. The authors of  \citep{kovacsvizvaririedeltoth} worked with 19 species, including four hypothetical transient species, made up of four atomic constituents (and charge); these are listed in Table \ref{tbl:KVRT-species}. The corresponding atomic matrix (see part I) is a $5\times 19$ matrix. To generate all combinatorially feasible elementary reactions, $2\cdot 19+{19 \choose 2} = 209$ systems need to be solved. The results give $1022$ combinatorially feasible elementary steps. Based on chemical evidence a large number of them were eliminated, leaving 673 elementary reactions.
\begin{table}[tb]
\center
\hrule
\(\mathbf{H_2C_2O_4}\) \hspace{1em}  \(\mathrm{HC_2O_4^-}\) \hspace{1em}  \(\mathrm{H}^+\) \hspace{1em}  \(\mathrm{C_2O_4^{2-}}\) \hspace{1em} \(\mathbf{Mn}^{2+}\) \hspace{1em} \(\mathrm{MnC_2O_4}\) \hspace{1em} \(\mathbf{MnO_4^-}\) \hspace{1em} \(\mathbf{MnO_2}\)  \hspace{1em}
\(\mathrm{Mn}^{3+}\) \hspace{1em}  \(\mathrm{CO}_2\) \hspace{1em} \(\mathrm{H_2O}\) \hspace{1em}  \(\mathrm{CO_2^-}\)
\([\mathrm{MnO_2,H_2C_2O_4}]\) \hspace{1em} \([\mathrm{Mn(C_2O_4)]^+}\) \hspace{1em} \([\mathbf{Mn(C_2O_4)_2}]^-\)
\([\mathrm{MnC_2O_4,MnO_4^-,H^+}]\) \hspace{1em} \([\mathrm{MnC_2O_4^{2+},MnO_3^-}]^+\) \hspace{1em} \([\mathrm{MnC_2O_4^{2+},MnO_3^-,H^+]^{2+}}\) \hspace{1em} \(\mathrm{[H^+,MnO_2,H_2C_2O_4]^+}\)
\hrule
\caption{Species of the permanganate/oxalic acid reaction according to \cite{kovacsvizvaririedeltoth}. Species with initially positive concentration are typeset bold.}
\label{tbl:KVRT-species}
\end{table}

The same computational results were reproduced in \citep{pappvizvari}, where further combinatorial analysis revealed that only 297 of the previously generated elementary reactions may be part of decompositions, and that 3 of the species cannot be generated during the reaction. This analysis takes into account which five species are present initially in the reaction (based on experiments), and is an adaptation of the indexing process by \cite{volpert}. It is perhaps worth mentioning that the result (in this specific example) is not sensitive to the assumption on the initial species; the same elementary steps and species remain if every non-complex species is assumed to be present initially in the vessel.

Decompositions of the overall reaction can be obtained via the solution of another linear Diophantine system of equations. With the above analysis the dimensions of this system are reduced to $16\times 297$. The generation of decompositions can be further accelerated by preprocessing: simple analysis, using linear programming, shows that every decomposition consists of at least 15 steps, and that two elementary steps take part in every decomposition, one of them with a coefficient of at least four. After this preprocessing, it is possible to find all decompositions of 15--17 steps, and thousands of more complex decompositions, from which chemically acceptable ones can be selected.

We refer the reader to the two papers cited above for details of these computations, and a review of a number of algorithms for the solution of linear Diophantine equations. All of these algorithms and the various methods of combinatorial analysis mentioned above are now built in to our program.

\section{Numerics}

\subsection{Stiffness}

One of the major problems when solving \ivp s describing chemical reactions is---as in general---stiffness:
the situation when the reaction proceeds along more than one time scales,
reaction rates of different steps are of different magnitude.
Then the codes meet the dilemma: if they choose a small step size, then the numerical solution will need much time, if they take a too long step size, then interesting phenomena occurring in the time evolution of the fast species will be lost. Starting with Gear \citep{gear} and Butcher \citep{butchercash} many methods fulfilling the requirements of the chemist have been elaborated and built in into program packages.
The following example show that such an insertion has been carried out especially successfully in the case of \Mma.

\begin{center}
\begin{figure}\label{robertsonfig1}
\includegraphics[height=5cm,width=9cm]{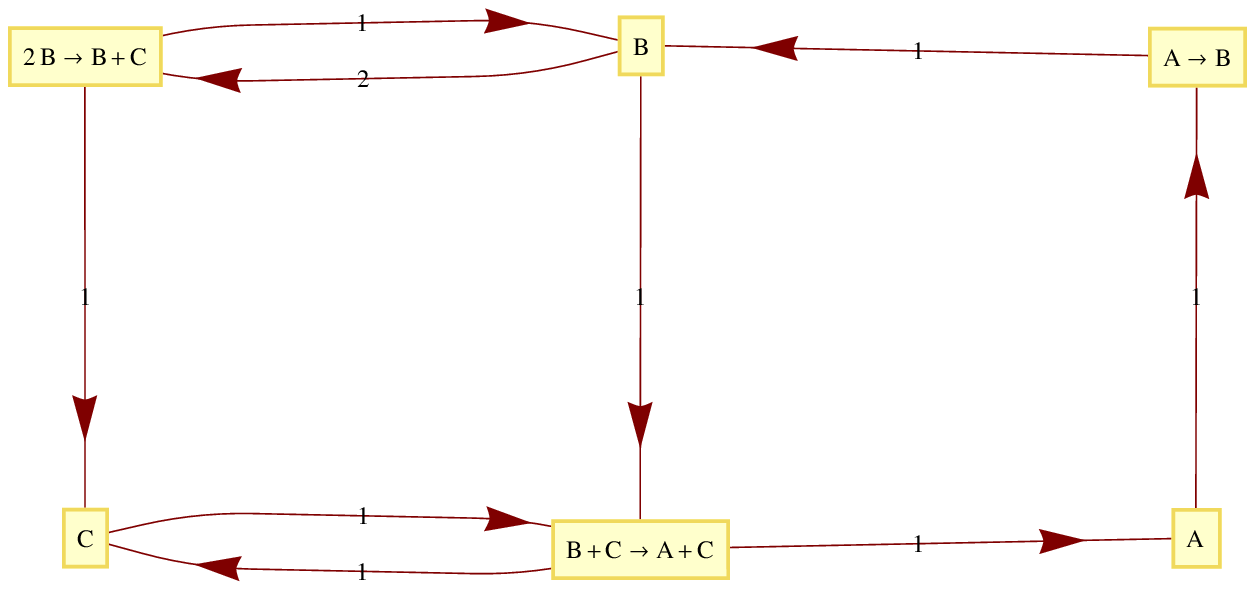}
\caption{The \Vg\ of the Robertson model}
\end{figure}
\end{center}

\begin{center}
\begin{figure}[h!]\label{robertsonfig31}
\includegraphics[width=8cm]{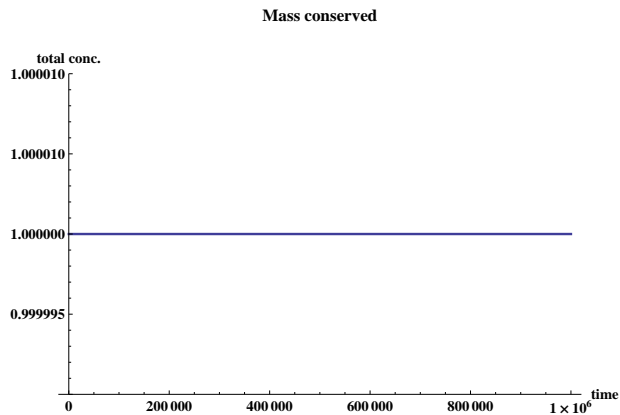}
\caption{The total mass in the Robertson model}
\end{figure}
\end{center}

The Robertson problem \citep{robertson}
\[
\irrevreak{\mathrm{A}}{0.04}{\mathrm{B}},\;
\irrevreak{2 \mathrm{B}}{3\times10^7}{\mathrm{B} + \mathrm{C}},\;
\irrevreak{\mathrm{B} + \mathrm{C}} {10^4}{ \mathrm{A} + \mathrm{C}}
\]
is a popular benchmark problem of stiff type, because there is a nine order of magnitude difference between the reaction rate coefficients.
To make matters worse, the solution is required on the time interval $[0,10^{11}].$
Figure \ref{robertsonfig1} shows its \Vg.
No matter what the values of the reaction rate coefficients are, the derivatives of the concentrations sum up to zero, meaning that total mass is conserved. One way to verify this fact is to execute the following:
\begin{center}
\textbf{\texttt{
Total[RightHandSide[{"Robertson"}, {k$_1$, k$_2$, k$_3$}]],
   }}
\end{center}
which gives 0, cf. \citep{yildirimbayramcons}.
The concentrations are calculated numerically as simply as usual (with our built-in function {\textbf{\texttt{Concentrations}}}):
\begin{center}
\textbf{\texttt{Concentrations[\{"Robertson"\},~\{0.04,~3$\times$10$^7$,~10$^4$\}, \{1, 0, 0\}, \{0,$10^{11}$\}]}}
\end{center}
Also the numerical method is good enough to conserve mass (Figure \ref{robertsonfig31}).
Let us see the individual concentrations of species A and C using the logarithmic time scale (Figure \ref{robertsonfig41}). But what about the species B? See Figure \ref{robertsonfig51}.

One of our colleagues, Dr. R. Horv\'ath, was so kind as to solve the same problem using MATLAB;
the same results were obtained within the same CPU time, but much less automatically.

Finally, we remark that a symbolic approach to unveil stiffness
when it does not manifest itself on the surface has been given by Prof.
Goldshtein at the conference, and see also \citep{bykovgoldshteinmaas}
\begin{center}
\begin{figure}[ht]\label{robertsonfig41}
\includegraphics[width=8cm]{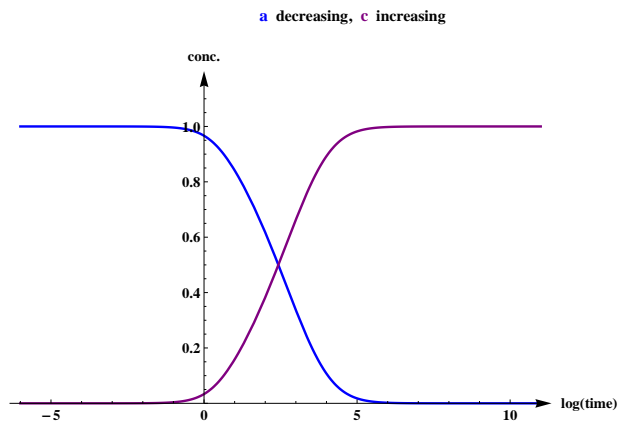}
\caption{Species A and C in the Robertson model}
\end{figure}
\vspace{3.5cm}
\begin{figure}[ht]\label{robertsonfig51}
\includegraphics[width=8cm]{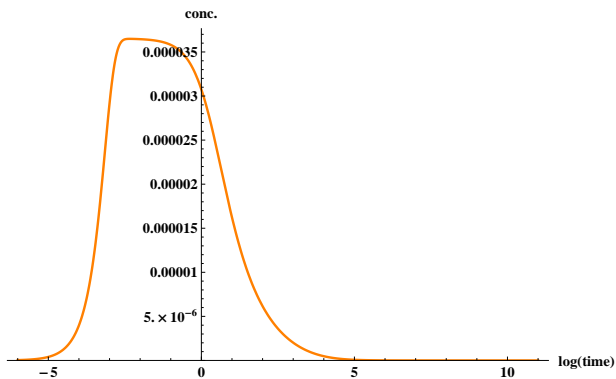}
\caption{Species B in the Robertson model}
\end{figure}
\end{center}

\subsection{Bifurcations}
A usual method to show that oscillatory solutions exist in the \ikde\ of a reaction is to apply
the well-known Andronov--Hopf bifurcation theorem.
Visualizing the effect of the change of the parameters
is really simple using the \Manipulate\ command of \Mma\,
see also the last example in the first part of our paper.
If you are interested how easily this can be done you might have a look at the demonstration by \cite{vardaitoth}.
In fact the method used there is based on particular ideas and has not as yet been built in into our program package handling general models.

\section{Symbolic methods}
Symbolic solutions are easier to find nowadays, still it is a hard problem even in the case of stationary points, see our previous paper.
The most interesting developments are expected in these fields, e.g.
how to find automatically the equivalence (of some kind) of different kinetic models
\citep{mendezfemat}, how to show the possibility or impossibility
of transforming one dynamics into another by (say, orthogonal) transformations
\citep{tothharsPhys} using perhaps the theory of algebraic invariants of differential equations
\citep{halmschlagerszenthetoth}, how to find a minimal model with given properties within a class of reactions
\citep{wilhelm,tothharsTCA}, given an \ikde\ how to find a reaction with a given structure
\citep{szederkenyi,szederkenyihangos}, etc.

\section{Stochastic simulation}
In reaction kinetics, it has always been obvious from the theoretical point of view that in the case of small systems (see e.g. \cite{aranyitoth}),
and systems operating around an unstable stationary state (e.g. models of chirality, see e.g. \cite{barabastothpalyi}) are better described
with the standard continuous time discrete state stochastic model of chemical reactions, see e. g.
\citep[Chaptr 5]{erditoth}. 
However, it turned out to be much more relevant only recently as with the development of analytical techniques it became possible to measure individual molecules. 
That is the reason why simulation of the stochastic model forms a relevant part of our program package.

Here we are considering only the usual continuous time discrete state model which is a Markovian jump process.
Denote by $\Xb(t)$ the number of species present in the system at time $t$.
Now, our investigated process is defined to evolve in time so that the law
\begin{equation}\label{def1}
P\left(\mbox{$R_r$ takes place in $[t,t+h)$}\,|\,\Xb(t)\right)=\kappa_r\left(\Xb(t)\right) h+\epsilon(h) h
\end{equation}
is fulfilled where $\epsilon(h)\to 0$ if $h\to0$ and $\kappa_j$'s  ($r=1,2,\ldots,R$) are a kind of \emph{combinatorial functions} of the vector of numbers of species similar to but slightly different from the product of powers (\emph{Kurtz kinetics}) involving also the reaction rate coefficients.
From this assumption one can easily determine the \emph{infinitesimal generator} of this Markov process from which 
the \emph{Kolmogorov forward} and \emph{backward equations} as well as the \emph{master equation} are obtained. However, our formulation is equivalent to the following explicit representation:
\begin{gather}\label{pois}
\Xb(t)=\Xb(0)+\sum_{r=1}^R Z_r\left(\int_0^t\kappa_r(\Xb(s))\,\dr s\right)\gammab(\cdot,r)
\end{gather}
where the $Z_r$'s are independent, unit-rate Poisson processes.
Note that there are many known stochastic simulation algorithms for the presented problem, but we intend to summarize only a few of them, namely: the \emph{direct method} (\cite{sipostotherdi1,sipostotherdi2}): let $\overline{\kappa}(\xb):=\sum_{r=1}^R\kappa_r(\xb)$ and suppose $T>0$ is given, then the initialization step is: $\xb(0):=\xb_0\in\RM$, $t:=0$\; then the algorithm works as follows:

\noindent
$
\While[t<T,\\
\gammab:=\RandomChoice[\{\frac{\kappa_r(\xb(t))}{\overline{\kappa}(\xb(t))}; r=1,2,\dots,R\}{\tt \to}\\
\hfill\{R_r; r=1,2,\dots,R\}];\\
s:=\RandomVariate[\ExponentialDistribution[\overline{\kappa}(\xb(t))]];\\
\xb(t+s):=\xb(t)+\gammab; t=t+s]\\
$

Notice that the idea of the direct method relies on the definition \eqref{def1}. Also there exists several variants (e.g. \emph{first reaction method}) and improvements of this early method (see e.g. \cite{gibsonbruck}).
In what follows we briefly present the \emph{approximation methods}, which have their roots in formula \eqref{pois}. At the heart of these kinds of methods is the \emph{leap condition}, that is choose $\tau$ to be small enough so that the change in the state in $[t,t+\tau)$ causes no relevant change for the $\kappa_r$'s.
Again, the initialization step is  $\xb(0):=\xb_0\in\RM$, $t:=0$\; then the algorithm works as follows:

\noindent
$
\While[t<T,\\
\text{Choose\ } \tau \text{\ so as to satisfy the leap condition;}\\
\Do[p_r:=\RandomVariate[\PoissonDistribution[\kappa_r(\xb(t))\tau],\\
\{r,1,R\}];\\
\xb(t+s):=\xb(t)+\sum_{r=1}^{R}P_r(\xb(t),\xb(t+\tau),p_r);\\
t=t+\tau]\\
$

\begin{center}
\begin{figure}[h!]\label{figbrus}
\includegraphics[width=8cm]{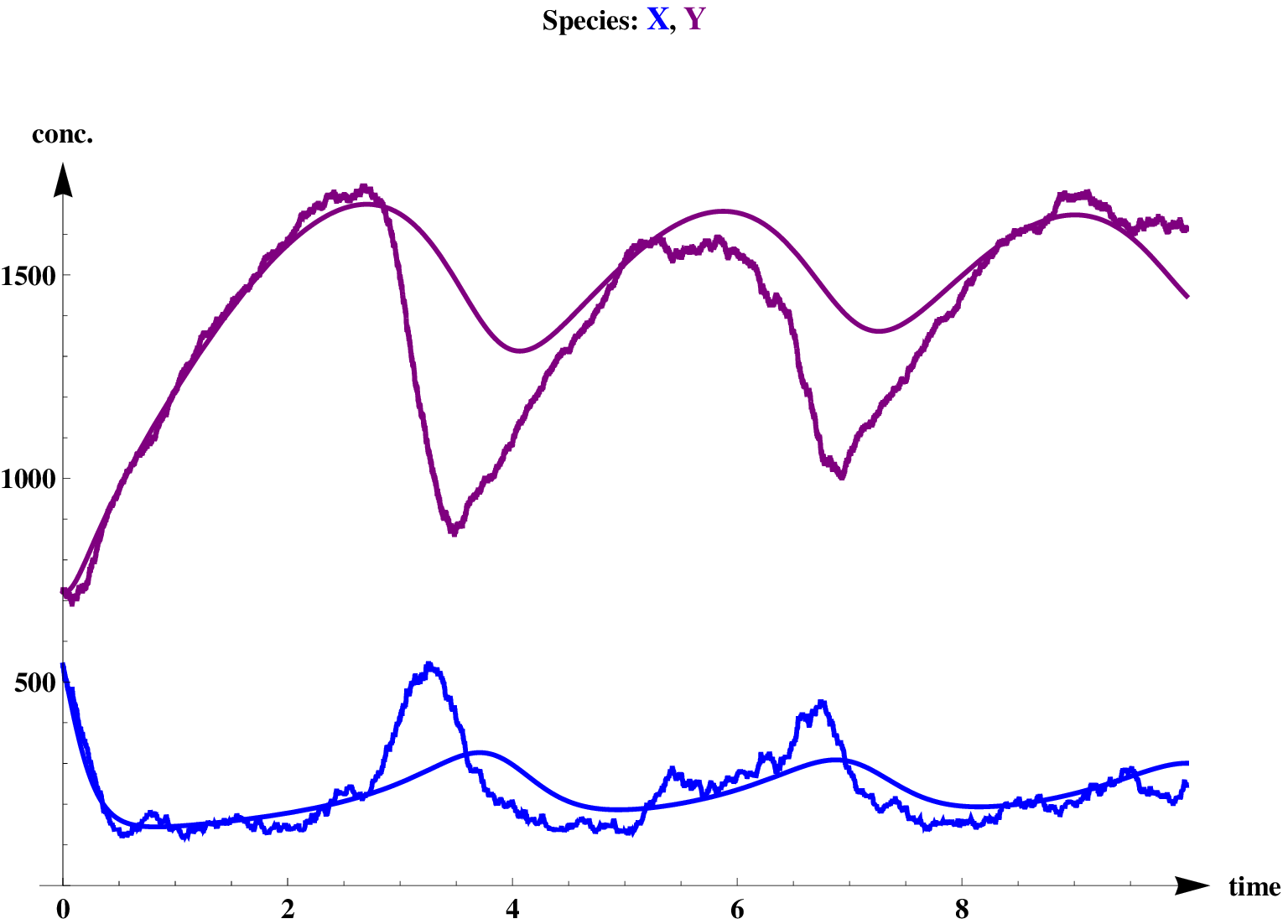}
\caption{Stochastic Brusselator model is considered with reaction rate coefficients: 
$k_1=1.92$, $k_2=5.76$, $k_3=5.6$, $k_4=4.8$ and initial conditions $x_0=500$, $y_0=720$. 
The ''noisy'' curve describes the evolution of the stochastic model, where the volume is $10^{-21}\mathrm{dm}^3$.}
\end{figure}
\end{center}

\begin{itemize}
\item 
If we take $P_r(\xb(t),\xb(t+\tau),p_r):=p_r,$ then we get the \emph{explicit $\tau$-leaping} method. 
\item
If $P_r(\xb(t),\xb(t+\tau),p_r):=p_r-\tau\kappa_r(\xb(t))+\kappa_r(\xb(t+\tau)),$ 
we are led to the \emph{implicit $\tau$-leaping} method.
\item
If we set  $P_r(\xb(t),\xb(t+\tau),p_r):=p_r-\frac{\tau}{2}\kappa_r(\xb(t))+\frac{\tau}{2}\kappa_r(\xb(t+\tau))$ then what we get is the \emph{trapezoidal $\tau$-leaping} method (\cite{cao,rathinam}).
\end{itemize}
However, these methods may also be considered as the stochastic analogues of certain numerical schemes applied for \ode s. Choosing the appropriate $\tau$ has turned out to be a non-trivial task. 
It is also a challenging problem to avoid negative population for $\xb(t)$ in the simulation process, see \citep{caopetzold}. 
Especially these last methods are used to handle \emph{stiff problems}, where typically implicit approaches give the most appropriate results.
These and several other algorithms have also been implemented into our program package. We also included conversion of units functions which provide an elegant way to compare the solution of the induced kinetic differential equation with the results of the stochastic simulation process which is considered to be taking place in a certain volume.
In Figure \ref{figbrus} we can see this comparison on the example of Brusselator model:
\[
\irrevreak{\mathrm{A}}{k_1}{\mathrm{X}},\irrevreak{\mathrm{X}}{k_2}{\mathrm{Y}},\irrevreak{2\mathrm{X} + \mathrm{Y}}{k_3}{3\mathrm{X}},\irrevreak{\mathrm{X}}{k_4}{\mathrm{P}}
\]
where A and P are external species.

\begin{center}
\begin{figure}[h!]\label{figauto}
\includegraphics[width=8cm]{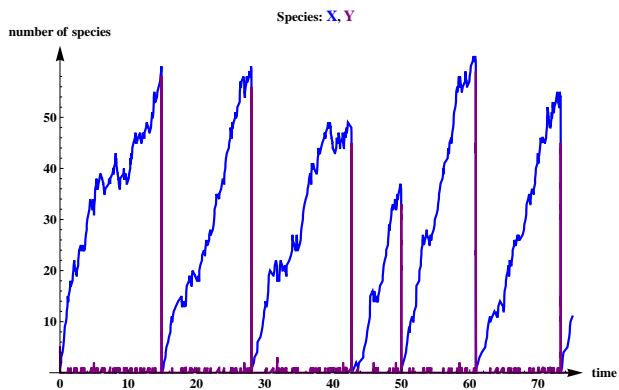}
\caption{Stochastic Autocatalator model with reaction rate coefficients: $k_1=110.2$, $k_2=0.094$, $k_3=0.011$, $k_4=90.34$, where the  initial conditions $x_0=15$, $y_0=80$ were considered.}
\end{figure}
\end{center}

In Figure \ref{figauto} one can see the behaviour of the Autocatalator model:
\[
\irrevreak{\mathrm{A}}{k_1}{\mathrm{X}},\irrevreak{\mathrm{X}}{k_2}{\mathrm{Y}},\irrevreak{\mathrm{X} + 2\mathrm{Y}}{k_3}{3\mathrm{Y}},\irrevreak{\mathrm{Y}}{k_4}{\mathrm{P}},
\]
where again A and P are external species.

\section{Estimation methods}
\subsection{Based on the deterministic model}
Estimating the reaction rate coefficients (or even the Arrhenius parameters $k_0, n$
and $A$ in the formula
\[
k(T)=k_0 T^n \exp(-A/(R T)),
\]
see \citep{nagyturanyi}) based on measurements is a kind of art. Here we only make a few remarks.

First of all, \Mma\ gives the possibility to obtain estimates (seemingly)
without any kind of iteration: it is capable using the numerical solution of a differential equation
in the same way as an explicitly given function.
Let us start from the deterministic model of the reaction
\[
\revreak{2\mathrm{X}\,}{k_1}{k_2}{\,\mathrm{X}},
\]
where $k_1=0.33$ and $k_2=0.72$ with the following initial concentration of X: \(x(0)=2.\)
Let us generate experimental data from the numerical solution and add a small amount of error.\\
\textbf{\texttt{
end = 7;
sol = First[
   x /.\\ NDSolve[\{x'[t] == -0.33 x[t]$^2$ + 0.72 x[t], \\
   \hfill x[0]~== 2\}, x, \{t, end\}]];\\
times = N[Range[0, 5 end]/5];\\
data = Transpose[\{times, sol[times] + RandomReal[.01, 5 end + 1]\}];\\
}}

\noindent Then, the \textbf{\texttt{model}} needed by the function \FindFit\\
\textbf{\texttt{FindFit[data, model[a, b][x],\\
\hfill \{\{a, 0.7\}, \{b, 0.2\}\}, x]}}\\
is defined in the following way.\\
\textbf{\texttt{
model[a\_,b\_] := (model[a, b] =
    First[x /.
      NDSolve[x'[t] == b x[t] - a x[t]$^2$, x[0] == 2\},
       x, \{t, end\}]])
       }}

\begin{center}
\begin{figure}[h!]\label{figfit}
\includegraphics[width=8cm]{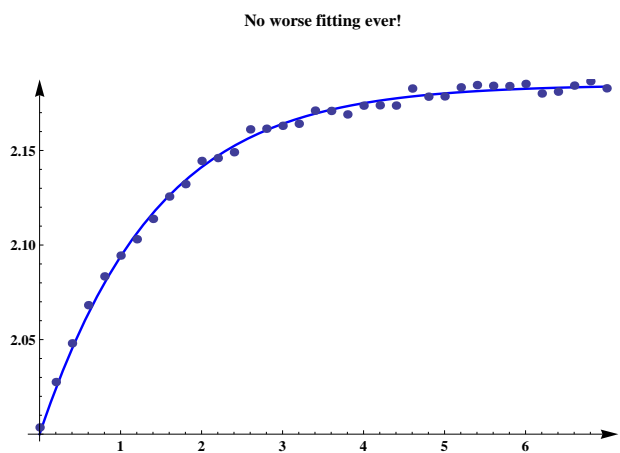}
\caption{Fitted model and data}
\end{figure}
\end{center}

As a result we get $(0.344129, 0.75174)$ instead of $(0.33,0.72),$
not a bad result, but certainly this is only a toy example.
The agreement between the ``measurements'' and the fitted model is quite good, see Figure \ref{figfit}.
A systematic investigation of similar (and further) estimation procedures for more complicated models
is in progress.
We also mention that a method similar to simulated annealing applied to the determination of Arrhenius parameters of reaction rate coefficients has also been elaborated \citep{nagyszaboturanyitoth}.

\subsection{Based on the stochastic model}
Based on the stochastic model one can obtain estimates in different ways.
If we are able to collect data on the individual changes of the molecules what is not so 
impossible in these days, then we can have maximum likelihood estimates. What is more even equilibrium fluctuations
might be enough to estimate the parameters, \citep[Subsection 5.8.7]{erditoth}.
Another approach is to use a kind of implicit linear regression: \cite{hangostoth}. 
Methods elaborated to estimate parameters of mass communication can also be adapted to the goals of reaction kinetics, see \cite{kovacs}.
\section{Visualization}
\begin{center}
\begin{figure}[h!]\label{figloc}
\centering
\includegraphics[height=8cm,width=7cm]{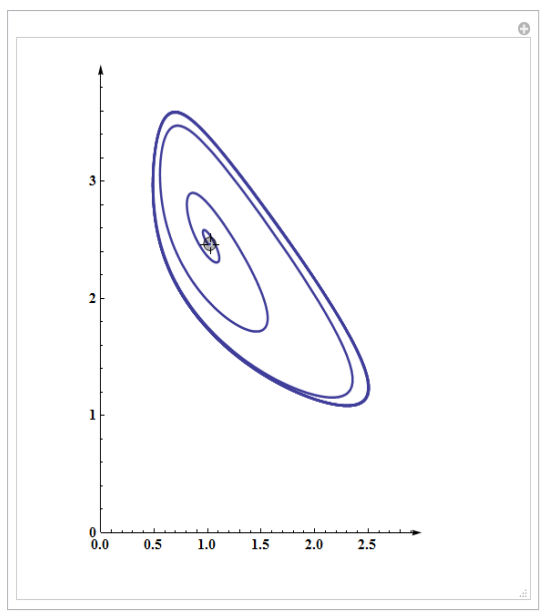}
\caption{The deterministic Brusselator model plotting in 2D using \Manipulate\ in which one can change the initial conditions in real time (using the ``locator'' object)}
\end{figure}
\end{center}

No external package is needed to provide camera ready figures for papers in any of the following forms:
eps, jpg, pdf, bmp, gif, wmf etc. (see a delicate example on Figure \ref{figloc}).

\section{Discussion, further plans}
We allow non mass action type kinetics even now,
therefore it does not seem to be a hard task to include
temperature dependence, so important in the case of modelling
e.g. atmospheric chemistry and combustion. Surely, reaction diffusion equations should also be treated in full generality including
symbolic treatment, as well. The reduction of the number of variables is also a very important topics, see e.g. \citep{tothlirabitztonmlin}.

Finally, please visit the website: \url{demonstrations.wolfram.com} and look for the demonstrations of Attila L\'{a}szl\'{o} Nagy and Judit V\'{a}rdai.

\subsection{Acknowledgments}
The authors thank the cooperation with dr. R. Horv\'ath.
\torol{
are very thankful for many remarks and corrigendums to Benedek Kov\'{a}cs, Ilona Nagy, Rudolf Csikja, \'{A}gota Busai and Tam\'{a}s Ladics.}
Many of the participants of MaCKiE 2011 contributed with incentive ideas. 
Further requirements, criticism and problems to be solved are wanted.

\bibliography{TothNagyPapp}%
\bibliographystyle{plainnat}
\end{document}